\documentclass{amsart}

\DeclareMathOperator{\core}{core}
\DeclareMathOperator{\rad}{rad}

\usepackage{amssymb}
\newtheorem{defn}{Definition}
\newtheorem{rem}{Remark}

\numberwithin{equation}{section}

\begin{document}

\title[Dirichlet Generating Functions]{Survey of Dirichlet Series of Multiplicative Arithmetic Functions}

\author{Richard J. Mathar}
\urladdr{http://www.strw.leidenuniv.nl/~mathar}
\email{mathar@strw.leidenuniv.nl}
\address{Leiden Observatory, Leiden University, P.O. Box 9513, 2300 RA Leiden, The Netherlands}

\subjclass[2010]{Primary 11K65, 11Y70; Secondary 30B50, 11M41}

\date{\today}
\keywords{Arithmetic Function, multiplicative, Dirichlet Generating Function}

\begin{abstract}
The manuscript reviews
Dirichlet Series of
important multiplicative
arithmetic functions. The aim is to represent these as products and ratios of
Riemann $\zeta$-functions, or, if that concise format is not found,
to provide the leading factors of the infinite product over $\zeta$-functions.
If rooted at the Dirichlet series for powers, for sums-of-divisors
and for Euler's totient,
the inheritance of multiplicativity through Dirichlet convolution or ordinary
multiplication of pairs of arithmetic functions generates most of the 
results.
\end{abstract}

\maketitle

\section{Scope}
\subsection{Definition}

Multiplicative functions are arithmetic functions $a(n)$---functions defined
for integer argument
$n\in \mathbb Z$---for which evaluation
commutes with multiplication for coprime arguments:
\begin{defn} (Multiplicative function $a$)
\begin{equation}
a(nm)= a(n)a(m) \quad \forall (n,m)=1.
\label{eq.Mdef}
\end{equation}
\end{defn}
The topic of the manuscript is the computation of the Dirichlet
series $\zeta_D$ of arithmetic functions of that kind for sufficiently
large real part of the argument $s$:
\begin{defn} (Dirichlet generating function $\zeta_D$)
\begin{equation}
a(n)\mapsto \zeta_D(s) \equiv \sum_{n\ge 1}\frac{a(n)}{n^s}.
\label{eq.zetaDdef}
\end{equation}
\end{defn}
As an immediate consequence of the definition,
the $\zeta_D$ of the product of $n^k$ times
a multiplicative function $a(n)$ is given by replacing $s\to s-k$ in the $\zeta_D$ of $a(n)$:
\begin{equation}
n^ka(n)\mapsto \zeta_D(s-k).
\label{eq.shift}
\end{equation}

\subsection{Properties.}
For ease of reference further down,
we summarize well-known features of
arithmetic functions.

A consequence of the definition
(\ref{eq.Mdef})
is that
the function is already entirely defined if specified for 
prime powers $p^e$, because all remaining values follow by
the prime power factorization of the arguments:
\begin{equation}
a(p_1^{e_1}p_2^{e_2}\cdots p_m^{e_m})= \prod_p a(p^e).
\end{equation}
The equation that explicitly specifies values of $a(p^e)$
will be called the \emph{master equation} of that sequence
in the sequel.

If $a$ is multiplicative, the Dirichlet series reduces to
a product over all primes $p$:
\begin{equation}
\zeta_D(s)
= \prod_p \left(1+ \frac{a(p)}{p^s}+\frac{a(p^2)}{p^{2s}}+\cdots\right)
= \prod_p \left(1+ \sum_{e\ge 1}\frac{a(p^e)}{p^{es}}\right)
= \prod_p \sum_{e\ge 0}\frac{a(p^e)}{p^{es}}
.
\label{eq.genZ}
\end{equation}
This sum over $e\ge 0$ at some fixed $p$
will
be called the \emph{Bell series} of that $a(n)$.
The product over the primes $p$ may be a finite product and/or ratio
of cyclotomic polynomials of some power of $p$;
then it is rephrased as a finite
product
of Riemann zeta-functions. In the general case, this expansion
will lead to an infinite product, and will be represented
in the followup chapters in the format
\begin{equation}
\zeta_D(s) = \prod_p
\prod_{i=1}^\infty(1-S_ip^{l_i-u_is})^{\gamma_i},
\end{equation}
where the vector of the $S_i$ contains sign factors
$\pm 1$, and the $l_i$, $u_i$ and $\gamma_i$
are integers.
The natural order of the factors is smallest $u_i$ first, and if these are the same
largest $l_i$ first. This order stresses which terms put the tightest
constraints on the region of convergence in the plane of complex $s$.
With Euler's formula, these constituents are equivalent
to Riemann $\zeta$-functions:
\begin{equation}
\prod_p (1-p^{l_i-u_is})^{\gamma_i} = \zeta^{-\gamma_i}(u_is-l_i)
;\quad
\prod_p (1+p^{l_i-u_is})^{\gamma_i}
=
\frac{
\zeta^{\gamma_i}(u_is-l_i)}{
\zeta^{\gamma_i}(2u_is-2l_i)}
.
\label{eq.Eulpm}
\end{equation}
(We write $\zeta_D$ for a generic Dirichlet series and $\zeta$
for the Riemann zeta function.)
Truncating the product expansion is
a tool of
numerical evaluation of the $\zeta_D(s)$.

The difference of this work to the Gould-Shonhiwa table 
\cite{GouldMJM20}
of the transformations $a(n)\mapsto \zeta_D(s)$
is that we will (i) cover Dirichlet series which require
this type of infinite
Euler products, and  (ii) will detail the six-digit A-numbers of
individual sequences in Sloane's Online
Encylopedia of Integer Sequences (OEIS) \cite{EIS};
due credit to individual's discovery of
many formulae that follow is stored in this database.

\subsection{Dirichlet convolution}
Dirichlet convolution is the construction of a new series
by summation over divisor and complementary divisor arguments of two arithmetic functions:
\begin{defn} (Dirichlet Convolution of $a$ and $b$)
\begin{equation}
(a\star b)(n)\equiv \sum_{d|n} a(d) b(n/d)
.
\label{eq.astarb}
\end{equation}
\end{defn}
The master equation of a multiplicative function derived via Dirichlet
convolution (\ref{eq.astarb}) is related to the master equations of
the factors  \cite{VaidyanTAMS33}:
\begin{equation}
(a\star b)(p^e) = \sum_{l=0}^e a(p^l)b(p^{e-l}).
\label{eq.astarbM}
\end{equation}

The ordinary (Hadamard) and the Dirichlet (convolution) product of two
multiplicative functions, the ordinary (Hadamard) ratio of two multiplicative
functions, and the Dirichlet inverse of a multiplicative function are multiplicative \cite{Apostol,BrownMathG89},
which creates multiplicative function by inheritance  from ``simpler'' multiplicative
functions.
As a special case, the $j-$th power of a multiplicative function is
multiplicative and the associated $a(p^e)$ is the
the $j$-th power of the one for the first power.

The Dirichlet series of the Dirichlet product is the (ordinary) product
of the Dirichlet series:
\begin{equation}
a\star b\mapsto \zeta_{D,a}(s)\zeta_{D,b}(s).
\label{eq.conv}
\end{equation}
As a corollary, the
Dirichlet series of the Dirichlet inverse defined by $(a^{(-1)}\star a)(n)\equiv \delta_{1,n}$
is given by the reciprocal Dirichlet series of $a$:
\begin{equation}
a^{(-1)} \mapsto 1/\zeta_D(s).
\label{eq.Dinv}
\end{equation}

\section{Classifications}
\subsection{Completely Multiplicative}
\subsubsection{Generic Properties}
Completely multiplicative functions
are a sub-species
of multiplicative functions which obey the equation (\ref{eq.Mdef}) for argument
pairs $n$ and $m$ irrespective of common divisors:
\begin{defn} (Completely Multiplicative function $a$)
\begin{equation}
a(nm)= a(n)a(m).
\label{eq.complMdef}
\end{equation}
\end{defn}
The well-known properties of completely multiplicative functions are
that the Dirichlet inverse can be written as a multiplication with the M\"obius function,
\begin{equation}
a^{(-1)}(n) = \mu(n)a(n),
\end{equation}
that the master equation allows interchange of exponentiation and evaluation,
\begin{equation}
a(p^e) = a^e(p),
\label{eq.Mpow}
\end{equation}
and that the sum over the exponents
$e$ in (\ref{eq.genZ}) is a geometric series \cite{Apostol},
\begin{equation}
\zeta_D = \prod_p\frac{1}{1-a(p)p^{-s}}.
\label{eq.eulComp}
\end{equation}

\subsubsection{Powers}
The Dirichlet series of powers is obvious from (\ref{eq.zetaDdef}):
\begin{equation}
n^k \mapsto \zeta(s-k).
\label{eq.Dpow}
\end{equation}
The most important example is the exponent $k=0$,
\begin{equation}
1 \mapsto \zeta(s).
\label{eq.Done}
\end{equation}
Its Dirichlet inverse  is the M\"obius function $\mu(n)$ (which is not
completely multiplicative, A008683), with Dirichlet generating function
\begin{equation}
\mu(n) \mapsto 1/\zeta(s)
\label{eq.Dmu}
\end{equation}
obtained combining (\ref{eq.Dinv}) and (\ref{eq.Done}).
The master equation of $\mu(n)$ is
\begin{equation}
a(p^e)=\left\{
\begin{array}{ll}
-1, & e=1;\\
0, & e>1.\\
\end{array}
\right.
\label{eq.Mmu}
\end{equation}
Squaring one obtains the Bell series of $\mu^2(n)$,
\begin{equation}
1+\sum_{e=1}\frac{1}{p^{es}}
= 1+p^{-s}
= \frac{1-p^{-2s}}{1-p^{-s}}
\end{equation}
which will be used further down as
\begin{equation}
\mu^2(n) \mapsto \zeta(s)/\zeta(2s).
\label{eq.Dmu2}
\end{equation}
\begin{rem}
These $n^k$ cover
$k=0$ in A000012,
$k=1$ in A000027 with inverse $n\mu(n)$ in A055615,
$k=2$ in A000290,
$k=3$ in A000578,
$k=4$ and $5$ in A000583 and A000584,
$k=6$--$9$ in A001014--A001017,
$k=10$--$12$ in A008454--A008456,
$k=13$--$25$ in A010801--A010813,
$k=26$ in A089081,
$k=27$-$30$ in A122968--A122971
\cite{EIS}.
\end{rem}

\subsubsection{Primes to constants}
If the master equation of a completely multiplicative function is
\begin{equation}
a(p) = c,
\end{equation}
the Dirichlet series is usually expanded into an Euler product
\begin{equation}
\zeta_D(s) = \prod_p\frac{1}{1-c/p^s} = \prod_{j\ge 1}(1-p^{-s})^{\gamma_j(c)}.
\label{eq.eulMc}
\end{equation}
for numerical efficienty,
such that the Dirichlet generating function becomes an (infinite) product
of the form $\prod_{j\ge 1} \zeta^{-\gamma_j(c)}(s)$ \cite{MoreeDM295}.
\begin{rem}
For $c=2$--$11$ the $\gamma_j(c)$ are
A001037,
A027376,
A027377,
A001692,
A032164,
A001693,
A027380,
A027381,
A032165,
and
A032166 in that order.
These exponents appear essentially
as $\gamma_{r,j}^{(C)}$ in
my expansions of Hardy-Littlewood constants
\cite[chapt.\ 7]{MatharArxiv0903}.
In numerical practise, Wynn's partial-product algorithm is
used to accelerate convergence of the Euler products \cite{WynnCJ14}.
\end{rem}
The cases of negative $c$ are mapped
via
\begin{equation}
\prod_p \frac{1}{1-cp^{-s}}
=
\prod_p \frac{1+cp^{-s}}{1-c^2p^{-2s}}
\end{equation}
to a division of two Hardy-Littlewood constants.

\begin{rem}
For a fixed integer $s$, one may factorize the polynomial
of $1/p$ over the reals numerically, to face a Weierstrass product representation
\begin{equation}
\zeta_D = \prod_p \prod_j (1+\frac{\beta_j}{p^t})
.
\end{equation}
The number of factors in the $j$-product is equivalent to the order of the
polynomial, and $\beta_j$ are essentially its roots.
Interchange of the two products rewrites $\zeta_D$ as a finite product
of prime zeta-functions of squarefree $k$-almost primes \cite{MatharArxiv0803}:
\begin{equation}
\zeta_D = \prod_j \left[1+\beta_j\sum_p \frac{1}{p^t} +\beta_j^2\sum_{p<q}\frac{1}{(pq)^t}+\cdots\right]
=
\prod_j \left[1+\sum_{k\ge 1}(-\beta_j)^k P_k^{(\mu)}(t)\right]
.
\end{equation}
\end{rem}
\begin{rem}
This covers A061142 and A165872 ($c=\pm 2$)
and
A165824 ($c=3$)
up to
A165871 ($c=50$) \cite{EIS}.
\end{rem}

\subsubsection{Liouville}
The Liouville function $\lambda(n)$ (A008836) is
the parity of the number $\Omega(n)$ of prime divisors of $n$. The master equation is \cite{SierpinskiMM42}
\begin{equation}
a(p) = -1,
\end{equation}
which evaluates by immediate application of (\ref{eq.eulComp}) to \cite{LehmanMC14}
\begin{equation}
\lambda(n) = (-1)^{\Omega(n)} \mapsto \zeta(2s)/\zeta(s),
\label{eq.liouv}
\end{equation}
the D-inverse of (\ref{eq.Dmu2}).
With (\ref{eq.conv}) follows that $\lambda \star 1$ is the characteristic function
of the squares \cite{SierpinskiMM42}:
\begin{equation}
\lambda \star 1 =  \epsilon_2(n)\mapsto \zeta(2s).
\end{equation}

\subsection{Persistently Multiplicative}
\subsubsection{Definition}
I call a multiplicative function $a$ \emph{persistently}
multiplicative if a product of coprime arguments leads
to a coprime product of the function,
\begin{defn} (Persistently multiplicative function $a$)
\begin{equation}
a(nm)=a(n)a(m)\quad \mathrm{and}\quad (a(n),a(m))=1\quad \forall (n,m)=1.
\end{equation}
\end{defn}
If a persistently multiplicative function $g$ is the inner function
of a compositorial product $a(n)=f(g(n))$, and if $f$ is multiplicative,
then $a$ is also multiplicative.

Persistently multiplicative are for example
those multiplicative
functions where master equations only modify the exponent
of the prime power through some function $E$, $a(p^e) = p^{E(e)}$.
The important subclass are the powers.
This also includes functions which remove all powers of some fixed prime $p_j$
from $n$, characterized by
\begin{equation}
\zeta_D = \zeta(s-1)(1-p_j^{1-s})/(1-p_j^{-s}).
\end{equation}
\begin{rem}
Examples with $p_j=2$, $3$ or $5$ are
A000265,
A038502,
and
A132739.
\end{rem}
Other persistently multiplicative functions permute
the prime bases $p$ in the master equation, for example replace primes
by their successors (A003961)
or swap with the adjacent prime (A061898).

\subsubsection{Squarefree core}
Persistently multiplicative are functions that reduce $n$ to its squarefree ($t=2$),
cubefree ($t=3$) etc
cores, where $E(e) = e\mod t$ is a modulo function which partitions $e$
into periodically modulated classes.
The function which reduces $n$ to the $t$-free core has the
Bell series
\begin{equation}
\sum_{r=0}^{t-1} \sum_{e=r,r+t,r+2t,\ldots} \frac{p^r}{p^{es}}
=
\sum_{r=0}^{t-1} \sum_{j=0}^\infty \frac{p^r}{p^{(r+jt)s}}
=
\frac{1-p^{t(1-s)}}{(1-p^{-ts})(1-p^{1-s})}
\label{eq.Bcore}
\end{equation}
and therefore with (\ref{eq.Eulpm}) the generating function
\begin{equation}
\core_t(n) \mapsto \zeta(ts)\zeta(s-1)/\zeta(ts-t).
\label{eq.Dtfreecore}
\end{equation}
\begin{rem}
This concerns sequences A007913 (squarefree), 
A050985 (cubefree) and A053165 (4-free).
\end{rem}

\subsubsection{Largest $t$-free Divisor}

The largest $t$-free number dividing $n$, $\rad_t(n)$, is complementary to
the functionality of the previous subsection.
The master equation admits exponents limited by $t$ and by the exponent in $n$:
\begin{equation}
a(p^e)= p^{\min(e,t-1)}.
\label{eq.Mradt}
\end{equation}
The Bell series is
\begin{multline}
\sum_{e=0}^{t-1} \frac{p^e}{p^{es}}
+\sum_{e\ge t} \frac{p^{t-1}}{p^{es}}
=
\frac{1-p^{-s}-p^{t(1-s)}+p^{t(1-s)-1} }{(1-p^{-s})(1-p^{1-s})}
\\
=
\frac{(1-p^{-s})\sum_{l=0}^{t-2}p^{(1-s)l}+p^{(1-s)(t-1)} }{1-p^{-s}}
.
\end{multline}
The denominator contributes $\zeta(s)$ to the Dirichlet series.
For $t=2$, the Euler expansion of the numerator starts:
\begin{multline}
\prod_p (1+p^{1-s}-p^{-s})
=
\prod_p
(1+p^{1-s})(1-p^{-s})(1+p^{1-2s})(1-p^{2-3s})
\\ \times
(1+p^{1-3s})
(1+p^{3-4s})(1-p^{2-4s})(1+p^{1-4s})(1-p^{4-5s})(1+p^{3-5s})^{2}
\\ \times
(1-p^{2-5s})^{2} (1+p^{1-5s})\cdots,\quad s>2
.
\label{eq.D2freed}
\end{multline}

\begin{rem}
The cases $t=2$--3 are shown in A007947--A007948,
the case $t=4$ in A058035. $1\star \rad_2$ is A191750.
\end{rem}

\subsubsection{Even-odd Splitting}
Persistently multiplicative are the functions that assign 1 to all odd
arguments and some other values to even arguments. A fundamental example
maps all even arguments to some constant $c$, which creates
an arithmetic sequence of period length 2:
\begin{equation}
a(p^e) =\left\{\begin{array}{ll}
c,& p=2;\\
1,& p>2;
\end{array} \right\}
\mapsto [1+(c-1)\cdot 2^{-s}]\zeta(s).
\end{equation}
In a variant, multiples of 4 could be assigned to some constant $c_1$, the other
even arguments to another constant $c_2$:
\begin{equation}
a(n) =\left\{\begin{array}{ll}1,& n\,\mathrm{odd} \\
c_1,& n\equiv 0\bmod 4\\
c_2,& n\equiv 2 \bmod 4 \\
\end{array} \right\}
\mapsto 
[1+(c_2-1)\cdot 2^{-s}+(c_1-c_2)\cdot 4^{-s}]\zeta(s)
.
\end{equation}
These periodic functions are additive overlays of
$L$-series \cite{MatharArxiv1008,BradleyRJ6}.
The computational strategy usually involves subtracting
the Riemann $\zeta$-function, expansion of the remaining $a(n)-1$
into a discrete Fourier series, and writing each
component as a Hilbert zeta-function.
\begin{rem}
This applies to
A109008 $(c_1=4,c_2=2)$, 
A010121 $(c_1=4,c_2=1)$,
A010123 $(c_1=6,c_2=2)$,
A010130 $(c_1=10,c_2=1)$,
A010131 $(c_1=10,c_2=2)$,
A010137 $(c_1=12,c_2=5)$,
A010146 $(c_1=14,c_2=6)$,
A112132 $(c_1=7,c_2=3)$,
A010127 $(c_1=8,c_2=3)$,
A089146 $(c_1=4,c_2=8)$,
or
A010132 $(c_1=10,c_2=4)$.
\end{rem}

\section{Core Classes}\label{sec.core}
\subsection{Characteristic Function of $t$-th powers}\label{sec.charf}
The characteristic function $\epsilon_t(n)$ of the $t$-th powers equals $1$ if 
the argument is a $t$-th power of some positive integer $b$, $0$ otherwise \cite{VaidyanTAMS33}.
The Dirichlet generating function (\ref{eq.zetaDdef}) collects $1/b^{ts}$ summing over all $b\ge 1$:
\begin{equation}
\epsilon_t(n) \mapsto \zeta(ts).
\label{eq.Dchart}
\end{equation}
The application of (\ref{eq.shift}) with  (\ref{eq.Dpow}) yields
\begin{equation}
n^{k/2}\epsilon_2(n) \mapsto \zeta(2s-k).
\label{eq.Dz2sk}
\end{equation}
\begin{rem}
$\sqrt n \epsilon_2(n)$ is A037213.
\end{rem}

The characteristic function of the numbers which are $t$-free (which cannot
be divided by a non-trivial $t$-th power) shall be denoted
$
\xi_t(n)
$.
The master equation puts a cap on the maximum power admitted in each factor:
\begin{equation}
a(p^e) = \left\{
\begin{array}{ll}
1, & e<t; \\
0, & e\ge t.
\end{array}
\right.
\end{equation}
The Bell series is
\begin{equation}
\sum_{e=0}^{t-1} 1/p^{es} = \frac{1-p^{-st}}{1-p^{-s}},
\end{equation}
therefore
\begin{equation}
\xi_t(n) \mapsto \zeta(s)/\zeta(st)
\label{eq.Dsqrfree}
\end{equation}
and \cite{VaidyanTAMS33}
\begin{equation}
\xi_t(n)\star \epsilon_t(n) = 1.
\end{equation}

\begin{rem}
The case $t=2$ comprises $\mu^2(n)$, the characteristic function of squarefree integers (A008966) \cite{SierpinskiMM42},
the D-inverse of (\ref{eq.liouv}).
The derived $n\xi_2(n)$
is represented by the absolute values of A055615.
\end{rem}

\subsection{Depleted $\zeta$-functions}
Characteristic functions
of numbers which are not multiples of some prime power $q^k$ are multiplicative
with
\begin{equation}
a(p^e)=\left\{
\begin{array}{ll}
1,& \mathrm{if}\quad p\neq q; \\
1,& \mathrm{if}\quad p=q,\quad e< k; \\
0,& \mathrm{if}\quad p=q,\quad e\ge k.
\end{array}
\right.
\end{equation}
The Bell series is $1/(1-p^{-s})$ for all $p\neq q$ and
$
\sum_{e=0}^{k-1}1/p^{es}
= (1-p^{-sk})(1-p^{-s})
$
for $p=q$. The merger of both is
\begin{equation}
\delta_{q^k\nmid n}
\mapsto
(1-q^{-sk})\zeta(s)
.
\end{equation}
\begin{rem}
Examples are $q^k=2^1$ in A000035,
$2^2$ in A166486,
$2^3$ in A168181,
$3^1$ in A011655 (multiplied by $n$ in A091684),
$3^2$ in A168182,
$5^1$ in A011558 (multiplied by $n$ in A091703),
$7^1$ in A109720,
$11^1$ in A145568,
or any principal Dirichlet character modulo some prime.
\end{rem}

\subsection{Greatest Common Divisors}
The greatest common divisor $(n,c)$
with respect to a constant $c$ is periodic $(n+c,c)=(n,c)$ \cite[\S 8.1]{Apostol}
and multiplicative. (Periodicity is revealed by the Euclidean
algorithm which starting from $n+c$ on one hand or $c$ on the other
yields the same quotients and remainders
already after the first step of the algorithm.)

Let $c=\prod_p p^{e_c}$ specify the prime exponents of the constant;
then the master equation is
\begin{equation}
a(p^e) = p^{\min(e,e_c)}
.
\label{eq.Mgcd}
\end{equation}
The Bell series is again an exercise in geometric series \cite[0.113]{GR}\cite{HsuFQ38},
\begin{equation}
\sum_{e=0}^{e_c}\frac{p^e}{p^{es}} + \sum_{e>e_c}\frac{p^{e_c}}{p^{es}}
=
\frac{ 1-p^{-s}+p^{e_c}p^{-(e_c+1)s}(1-p) }
{(1-p^{-s})(1-p^{1-s})}
.
\end{equation}
The product over all primes, the Dirichlet series, is the
Riemann $\zeta$-function multiplied by a product of rational polynomials
over the primes with non-vanishing $e_c$:
\begin{multline}
(n,c)\mapsto
\zeta(s)
\prod_{e_c> 0}
\frac{ 1-p^{-s}+p^{e_c}p^{-(e_c+1)s}(1-p) }
{1-p^{1-s}}
\\
=
\zeta(s)
\prod_{e_c> 0}
\left(
1+(p-1)\sum_{l=0}^{e_c-1}p^{l(1-s)-s}
\right)
.
\label{eq.Dgcd}
\end{multline}
\begin{rem}
The reference sequences are A109007--A109015 for $c=$3--12 in the OEIS \cite{EIS},
with the exception of $c=6$ which is A089128.
\end{rem}

\subsection{Least Common Multiples}
The least common multiple $[n,c]$ of $n$ and a constant $c$ is constructed
with the master equation $a(p^e) = p^{\max(e,e_c)}$ but is not multiplicative in general.
With $(n,c)[n,c]=nc$ and multiplicativity of $(n,c)$, the divided
$[n,c]/c=n/(n,c)$ serves as a multiplicative substitute. The master equation of $[n,c]/c$ is
\begin{equation}
a(p^e)
= p^{ \max(e,e_c)}/p^{e_c} 
= p^{ \max(e-e_c,0)}.
\label{eq.Mlcm}
\end{equation}
The Bell series is
\begin{equation}
\sum_{e=0}^{e_c}\frac{1}{p^{es}}
+\sum_{e>e_c}\frac{p^{e-e_c}}{p^{es}}
=\frac{ 1-p^{1-s}+p^{-s(1+e_c)}(p-1)}{(1-p^{-s})(1-p^{1-s})}
.
\end{equation}
The analog of (\ref{eq.Dgcd}) becomes
\begin{equation}
[n,c]/c
\mapsto \zeta(s-1)
\prod_{e_c>0}
\left(
1+(1-p) \sum_{l=0}^{e_c-1}p^{-(l+1)s}
\right)
.
\end{equation}
\begin{rem}
This refers for $c=2$--$20$ to A026741,
A051176, A060819, A060791, A060789, A106608--A106612, A051724, and A106614--A106621.
\end{rem}

\subsection{Sigma: Sum of Divisors}
\subsubsection{Base Sequence}
The divisors of some number $n$
\begin{equation}
n = \prod _p p_i^{e_i}
\end{equation}
are of the form $d=\prod_p p_i^{m_i}$ with $0\le m_i\le e_i$.
The sum of the $k$-th power of divisors is
\begin{equation}
\sigma_k(n) = (1+p_1^k+p_1^{2k}+\cdots p_1^{e_1k})
(1+p_2^k+p_2^{2k}+\cdots p_2^{e_2k})\cdots
\end{equation}
which is a product of geometric sums \cite[p.\ 239]{HardyWright}:
\begin{equation}
\sigma_k(p^e) =
\left\{
\begin{array}{ll}
\frac{p^{k(e+1)}-1}{p^k-1},&\quad k>0; \\
e+1,& \quad k=0. \\
\end{array}
\right.
\label{eq.Mzeta}
\end{equation}
Inserted into (\ref{eq.genZ}) provides the Dirichlet series
\begin{equation}
\zeta_D =
\prod_p \left(\sum_{e\ge 0}\frac{p^{k(e+1)}-1}{p^k-1}\cdot \frac{1}{p^{es}}\right),
\quad k>0
,
\end{equation}
and the geometric series is summarized as \cite[(1.3.1)]{Titchmarsh}\cite[p.\ 293]{Dickson}
\begin{equation}
\sigma_k(n) \mapsto
\zeta(s)\zeta(s-k),\quad  k\ge 0
.
\label{eq.Dsigmak}
\end{equation}
In view of (\ref{eq.astarb}) and (\ref{eq.Dpow}) this shows
\begin{equation}
\sigma_k(n) = n^k \star 1.
\end{equation}

\begin{rem}
This covers
$k=0$, $1\star 1=\sigma_0(n)$, in A000005 with D-inverse A007427,
$k=1$ in A000203 with D-inverse A046692,
$k=2$ in A001157 with D-inverse A053822,
$k=3$ in A001158 with D-inverse A053825,
$k=4$ in A001159 with D-inverse A053826,
$k=5$ in A001160 with D-inverse A178448,
and $k=6$--$24$ in A013954--A013972.
\end{rem}
The sum over the inverse $k$-th powers deals with negative indices
of the $\sigma$-function. By inspection of the
complementary divisors $n/d$ for each $d$ this is
\begin{equation}
\sigma_{-k}(n) = \sum_{d|n} \frac{1}{d^k} = \frac{\sigma_k(n)}{n^k}.
\end{equation}
Applying the shift-theorem (\ref{eq.shift}) demonstrates
that (\ref{eq.Dsigmak}) is also valid in the range $k<0$.

\subsubsection{Convolutions}
With (\ref{eq.Dsigmak}) we derive for example
$\sigma_1(n)\star 1\mapsto \zeta^2(s)\zeta(s-1)$
(A007429),
$\sigma_0(n)\star \sigma_1(n)\mapsto \zeta^3(s)\zeta(s-1)$ (A007430),
$\sigma_2(n)\star 1\mapsto \zeta^2(s)\zeta(s-2)$
(A007433)
or
$\sigma_1(n)\star \sigma_1(n)\mapsto \zeta^2(s)\zeta^2(s-1)$ (A034761).

\subsection{Sums of Divisors which are $t$-th Powers}
The sum over all divisors of $n$ which are perfect $t-$th powers
is
\begin{equation}
a(n) = \sum_{d|n} d \epsilon_t(d) = n\epsilon_t(n) \star 1
\mapsto \zeta(s)\zeta(ts-t)
\label{eq.Ddivt}
\end{equation}
using the notation of the characteristic function $\epsilon_t$
(Section \ref{sec.charf}).
\begin{proof}
The Dirichlet generating function in (\ref{eq.Ddivt}) is derived
(i) either
by summing the Bell series and noting that the denominators of the
intermediate result are cyclotomic polynomials of $p^{-s}$ which allows
to express the Euler product as a finite product of $\zeta$-functions, or
(ii) more quickly starting from the
generating function
(\ref{eq.Dchart}) of $\epsilon_t(n)$, using the shift theorem
(\ref{eq.shift}) to produce the generating function for $n\epsilon_t(n)$,
\begin{equation}
n\epsilon_t(n)\mapsto \zeta(ts-t)
\label{eq.Dcharn}
\end{equation}
and exploiting the convolution with 1 via (\ref{eq.conv}) and (\ref{eq.Done}).
\end{proof}
\begin{rem}
The examples are $t=2$, the sums of the square divisors (A035316),
and $t=3$, the sum of the cube divisors (A113061).
\end{rem}
The 
master equation is
\begin{equation}
a(p^e) = \sum_{l= 0}^{\lfloor e/t\rfloor} p^{lt}
=\frac{p^{t(1+\lfloor e/t\rfloor)}-1}{p^t-1}
\label{eq.Msigtp}
\end{equation}
which can be made more explicit by writing
this down for each remainder of $e\mod t$ in the style of (\ref{eq.Bcore}).

The largest $t$-th power dividing $n=\prod_p p^e$ may be written
as $\max_{b^t\mid n}$.
For each prime basis $p$ it selects the maximum exponent $e$ which is a multiple of $t$.
This reduces the sum (\ref{eq.Msigtp}) over all multiples to its largest term:
\begin{equation}
a(p^e) = p^{t\lfloor e/t\rfloor}
\label{eq.Mltp}
\end{equation}
Substituting $e=kt+r$ in the Bell series yields
\begin{equation}
1+\sum_{e\ge 1}\frac{p^{e-r}}{p^{es}}
=
\sum_{r=0}^{t-1} \sum_{k\ge 0}\frac{p^{kt}}{p^{(kt+r)s}}
=
\sum_{r=0}^{t-1} p^{-rs}\frac{1}{1-p^{t(1-s)}}
=
\frac{1-p^{-st}}{(1-p^{-s})(1-p^{t(1-s)})}
.
\end{equation}
The product over all primes is
\begin{equation}
\max_{b^t\mid n} \mapsto  \zeta(s)\zeta(ts-t)/\zeta(st).
\end{equation}
Multiplying this $\zeta$-product by (\ref{eq.Dtfreecore}) shows
in conjunction with (\ref{eq.conv}) and (\ref{eq.Dsigmak}) that
\begin{equation}
\max_{b^t\mid n} \star \core_t(n) = \sigma_1(n).
\end{equation}
\begin{rem}
Examples are $t=2$, the largest square dividing $n$ (A008833),
$t=3$, the largest cube dividing $n$ (A008834),
or $t=4$, the largest 4th power dividing $n$ (A008835).
\end{rem}
One can also split the product in view of (\ref{eq.Dcharn}) and (\ref{eq.Dsqrfree}),
\begin{equation}
\max_{b^t\mid n} = \xi_t(n) \star (n\epsilon_t(n)).
\end{equation}

A similar function is the $t$-th root of the largest $t$-th power dividing
$n$,
\begin{equation}
a(p^e) = p^{\lfloor e/t\rfloor}
\label{eq.Metfloor}
,
\end{equation}
pulling the $t$-th root out of (\ref{eq.Mltp}).
Bell and Dirichlet series are
\begin{equation}
\sum_{e\ge 0} \frac{p^{\lfloor e/t\rfloor}}{p^{es}}
=
\sum_{r=0}^{t-1}\sum_{k\ge 0} \frac{p^k}{p^{(kt+r)s}}
=
\frac{1-p^{-st}}{(1-p^{-s})(1-p^{1-st})}
\mapsto
\zeta(st-1)\zeta(s)/\zeta(st).
\label{eq.Dmaxbtn}
\end{equation}
\begin{rem}
This theory applies to A000118 ($t=2$),
A053150 ($t=3$) and
A053164 ($t=4$).
\end{rem}

\subsection{Sum of $t$-free Divisors}
The sum of the $k$-th powers of $t$-free divisors of $n$ is
---in the notation of section \ref{sec.charf}---
\begin{equation}
\sum_{d|n} d^k\xi_t(d) = (n^k\xi_t(n))\star 1 \mapsto \zeta(s)\zeta(s-k)/\zeta(ts-tk),\quad k\ge 0.
\label{eq.Dsigtfree}
\end{equation}
This Dirichlet series follows applying (\ref{eq.shift}) to (\ref{eq.Dsqrfree})
and then (\ref{eq.Done}) and (\ref{eq.conv}).
\begin{rem}
The count of the squarefree divisors is A034444 with D-inverse in A158522;
the count of the cubefree divisors is A073184.
The sum of squarefree divisors (A048250)
has the master equation
\begin{equation}
a(p^e)= p+1.
\end{equation}
Multiplication by $n$ generates A181797.
The sum of the cubefree divisors is A073185.
\end{rem}

The count of the $t$-full divisors has the master equation \cite{SuryanAMS34}
\begin{equation}
a(p^e) = \max(1,e-t+2),
\label{eq.Msig0tfull}
\end{equation}
assuming $1$ is included in the set of $t$-full numbers.
Compared to the full count of divisors, this eliminates
contributions of the powers $p^1$, $p^2$,\ldots ,$p^{t-1}$
from the prime factorization of the divisors. The Bell series is
\begin{equation}
\sum_{e=0}^{t-1}\frac{1}{p^{es}}+\sum_{e\ge t}\frac{e-t+2}{p^{es}}
=
\frac{p^{-st}-p^{-s}+1}{(1-p^{-s})^2}
.
\end{equation}
For $t=2$, the numerator polynomial is the cyclotomic polynomial $\Phi_6(p^{-s})$,
and expansion of numerator and denominator with $1+p^{-s}$ yields
\begin{equation}
\sum_{d:\xi_2(d)=0} 1 \mapsto \zeta(s)\zeta(2s)\zeta(3s)/\zeta(6s).
\end{equation}
For $t=3$ and $t=4$ the infinite Euler products start as
\begin{equation}
\sum_{d:\xi_3(d)=0} 1 \mapsto \zeta(s)
\prod_p (1+p^{-3s})(1+p^{-4s})(1+p^{-5s})(1+p^{-6s})(1-p^{-9s})\cdots, \, s>1,
\end{equation}
and
\begin{equation}
\sum_{d:\xi_4(d)=0} 1 \mapsto \zeta(s)
\prod_p 
(1+p^{-4s})(1+p^{-5s})(1+p^{-6s})(1+p^{-7s})(1+p^{-8s})(1-p^{-11s})\cdots,\, s>1
.
\end{equation}
\begin{rem}
$t=2$ is A005361. $t=3$ is A190867.
\end{rem}

\subsection{Sigma of powers}
$\sigma_k(n^2)$ is an arithmetic function with master equation obtained by the
substitution $e\to 2e$ in (\ref{eq.Mzeta}):
\begin{equation}
a(p^e) = \frac{p^{k(2e+1)}-1}{p^k-1},\quad k>0;\quad a(p^e) =2e+1,\quad k=0.
\end{equation}
The Bell series is
\begin{equation}
\sum_{e\ge 0}\frac{p^{k(2e+1)}-1}{(p^k-1)p^{es}} = \frac{1+p^{k-s}}{(1-p^{2k-s})(1-p^{-s})},\quad k\ge 0,
\end{equation}
which induces
\begin{equation}
\sigma_k(n^2)\mapsto \zeta(s)\zeta(s-k)\zeta(s-2k)/\zeta(2s-2k).
\label{eq.Dsigkn2}
\end{equation}
If the right hand side is interpreted as the product of $\zeta(s)\zeta(s-2k)$ and
$\zeta(s-k)/\zeta(2s-2k)$, equations
(\ref{eq.Dsqrfree})
and
(\ref{eq.Dsigmak})
demonstrate
\begin{equation}
\sigma_k(n^2) = \sigma_{2k}(n) \star (n^k \xi_2(n) ).
\end{equation}
An alternative interpretation as a product of $\zeta(s-2k)$ and $\zeta(s)\zeta(s-k)/\zeta(2s-2k)$
shows with (\ref{eq.Dsigtfree})
\begin{equation}
\sigma_k(n^2) = [\sum_{d\mid n} d^k\xi_2(d)] \star n^{2k}.
\end{equation}
\begin{rem}
The case  $\sigma_0(n^2)$ in A048691 is documented
by Titchmarsh \cite[(1.2.9)]{Titchmarsh}
with
$\sigma_0(n^2)\mapsto  \zeta^3(s)/\zeta(2s)$.
$\sigma_1(n^2)$ is A065764,
and
$\sigma_2(n^2)$ is A065827.
\end{rem}

Moving on to higher powers in the argument, subsampled sums of divisors,
we first meet $\sigma_0(n^t)$
with Bell series
\begin{equation}
1+\sum_{e\ge 1}\frac{te+1}{p^{es}} = \frac{1+(t-1)p^{-s}}{(1-p^{-s})^2}
.
\label{eq.Bsigma0t}
\end{equation}
The denominator contributes a factor
$\zeta^2(s)$ to the Dirichlet series, and the numerator
is covered by division through the associated term of (\ref{eq.eulMc}).

The master equation of $\sigma_1(n^t)$
replaces $e$ by $et$ in (\ref{eq.Mzeta}),
\begin{equation}
a(p^e) = \frac{p^{et+1}-1}{p-1},
\end{equation}
which generates a Bell series
\begin{equation}
1+\sum_{e\ge 1}
\frac{p^{te+1}-1}{(p-1)p^{es}}
=
\frac{1+p^{1-s}\sum_{l=0}^{t-2}p^l}{(1-p^{-s})(1-p^{t-s})}
.
\label{eq.Bsig1nt}
\end{equation}
At $t=3$, the Euler expansion starts
\begin{multline}
\sigma_1(n^3)\mapsto  \zeta(s)\zeta(s-3)
\prod_p
(1+p^{2-s})(1+p^{1-s})(1-p^{3-2s})(1+p^{5-3s})
\\ \times
(1+p^{4-3s})
(1-p^{7-4s})(1-p^{6-4s})(1-p^{5-4s})(1+p^{9-5s})
\\ \times
(1+p^{8-5s})^{2}
(1+p^{7-5s})^{2} (1+p^{6-5s})\cdots
, \quad s>4,
\label{eq.Dsig1n3}
\end{multline}
for example. At $t=4$ it is
\begin{multline}
\sigma_1(n^4)\mapsto  \zeta(s)\zeta(s-4)
\prod_p
(1+p^{3-s})(1+p^{2-s})(1+p^{1-s})(1-p^{5-2s})
\\ \times
(1-p^{4-2s})
(1-p^{3-2s})(1+p^{8-3s})(1+p^{7-3s})^{2} (1+p^{6-3s})^{2}
\\ \times
(1+p^{5-3s})^{2}
(1+p^{4-3s})(1-p^{11-4s})(1-p^{10-4s})^{2} (1-p^{9-4s})^{4}
\\ \times
(1-p^{8-4s})^{4}
(1-p^{7-4s})^{4} (1-p^{6-4s})^{2} (1-p^{5-4s})\cdots, \quad s>5
.
\label{eq.Dsig1n4}
\end{multline}

\begin{rem}
Templates of these sequences are $\sigma_0(n^3)$ is A048785,
$\sigma_1(n^3)$ in A175926.
$\sigma_0(n^2)\star 1$ is A035116.
$\sigma_0(n^3)\star 1$ is A061391.
\end{rem}

\subsection{Sum of Gcd or Lcm}

Following (\ref{eq.Mgcd}),
the gcd of a divisor $d$ and its complementary divisor $n/d$ contributes
with a factor $(p^m,p^{e-m})^t=p^{t\cdot \min(m,e-m)}$ to $\sum_{d|n} (d,n/d)^t$.
Summing over $m$ from $0$ to $e$ yields
the master equation
\begin{equation}
a(p^e) = \left\{
\begin{array}{ll}
p^{et/2}+2\sum_{m=0}^{e/2-1} p^{tm}=[(p^t+1)p^{et/2}-2]/(p^t-1),& e\,\mathrm{even};\\
2\sum_{m=0}^{(e-1)/2} p^{tm}= 2[p^{t(e+1)/2}-1]/(p^t-1),& e\,\mathrm{odd}.\\
\end{array}
\right.
\end{equation}
The Bell series is
\begin{equation}
\sum_{e=0,2,4,\ldots} \frac{(p^t+1)p^{et/2}-2}{(p^t-1)p^{es}}
+\sum_{e=1,3,5,\ldots} 2\frac{p^{t(e+1)/2}-1}{(p^t-1)p^{es}}
=
\frac{1+p^{-s}}{(1-p^{-s})(1-p^{t-2s})}
,
\end{equation}
which reveals
\begin{equation}
\sum_{d|n} (d,n/d)^t
\mapsto
\zeta^2(s)\zeta(2s-t)/\zeta(2s)
\label{eq.Dsigmagcd}
.
\end{equation}

The associated analysis for the lcm starts from (\ref{eq.Mlcm}).
The prime $p$ and exponent $m$ of the divisor $d$ contribute to
$\sum_d [d,n/d]^t$ with a term $[p^m,p^{e-m}]^t=
p^{t\cdot \max(m,e-m)}$. The master equation splits again into two
cases depending on whether a middle term at $m=e/2$ is present or not:
\begin{equation}
a(p^e) = \left\{
\begin{array}{ll}
p^{et/2}+2\sum_{m=0}^{e/2-1} p^{t(e-m)}=[2p^{t(1+e)}-(p^t+1)p^{et/2}]/(p^t-1),& e\,\mathrm{even};\\
2\sum_{m=0}^{(e-1)/2} p^{t(e-m)}= 2e^{t(e+1)/2}[p^{t(e+1)/2}-1]/(p^t-1),& e\,\mathrm{odd}.
\end{array}
\right.
\end{equation}
The Bell series factorizes in the $\zeta$-basis
similar to (\ref{eq.Dsigmagcd}):
\begin{equation}
\sum_{d\mid n} [d,n/d]^t \mapsto
\zeta^2(s-t)\zeta(2s-t)/\zeta(2s-2t)
.
\end{equation}

\begin{rem}
$\sum_d (d,n/d)^t$ is A055155 for $t=1$ and A068976 for $t=2$.
$\sum_d [d,n/d]$ is A057670.
\end{rem}

\subsection{Sigma powers}
\nopagebreak
\subsubsection{Ordinary Products}

The $t$-th power of (\ref{eq.Mzeta}) is
\begin{equation}
\sigma_k^t(p^e) = 
\frac{(p^{k(e+1)}-1)^t}{(p^k-1)^t},\quad e\ge 0 .
\end{equation}
The binomial expansion of the associated Bell series is
\begin{multline}
\sum_{e\ge 0}\frac{\sigma_k^t(p^e)}{p^{es}}
=
\frac{1}{(p^k-1)^t}\sum_{e\ge 0}\sum_{t'=0}^t (-)^{t-t'}\binom{t}{t'}
\frac{ p^{t'k(e+1)}}{p^{es}}
\\
=
\frac{1}{(p^k-1)^t}\sum_{t'=0}^t (-)^{t-t'}\binom{t}{t'}\frac{p^{t'k}}{1-p^{(t'k-s)}}
\label{eq.Bsigmakt}
.
\end{multline}
For the squares of $\sigma$, $t=2$,
\begin{equation}
\sum_{e\ge 0}\frac{\sigma_k^2(p^e)}{p^{es}}
=
\frac{1-p^{2k-2s}}{(1-p^{-s})(1-p^{k-s})^2(1-p^{2k-s})}
\end{equation}
produces the Dirichlet series
\begin{equation}
\sigma_k^2(n)\mapsto 
\frac{\zeta(s)\zeta^2(s-k)\zeta(s-2k)}{\zeta(2s-2k)}
.
\label{eq.Dsigmak2}
\end{equation}
Because this equals (\ref{eq.Dsigkn2}) multiplied by $\zeta(s-k)$, we find
with (\ref{eq.shift}) and (\ref{eq.astarb}):
\begin{equation}
\sigma_k^2(n) = \sigma_k(n^2)\star n^k.
\end{equation}
\begin{rem}
These considerations cover A035116 with the Dirichlet series
\cite[(1.2.10)]{Titchmarsh}
\begin{equation}
\sigma_0^2(n) = \sigma_0(n^2) \star 1 = \sigma_0(n) \star \sigma_0^{*}(n) \mapsto  \zeta^4(s)/\zeta(2s),
\label{eq.Dsigma0sq}
\end{equation}
where $\sigma_0^{*}(n)$ is the number of unitary divisors of $n$ (A034444).
They also cover
$\sigma_1^2(n)$ in A072861.
\end{rem}
For larger $t$, the denominators of (\ref{eq.Bsigmakt})
contribute $\prod_{t'=0}^t \zeta(s-t'k)$
to the Dirichlet series (represented for $k=1$ and various $t$ by A001001 and A038991--A038999),
but the numerators do not factor as nicely.
The examples are
\begin{equation}
\sum_{e\ge 0} \frac{\sigma_k^3(p^e)}{p^{es}}
=
\frac{p^{3k-2s}+2p^{2k-s}+2p^{k-s}+1}
{\prod_{t'=0}^3(1-p^{t'k-s})}
\label{eq.Bsigmak3}
\end{equation}
or
\begin{equation}
\sum_{e\ge 0} \frac{\sigma_k^4(p^e)}{p^{es}}
=
\frac{p^{6k-3s}+(3p^{2k}+5p^{k}+3)p^{3k-2s}+(3p^{2k}+5p^{k}+3)p^{k-s}+1}
{\prod_{t'=0}^4(1-p^{t'k-s})}
.
\end{equation}

The Euler product expansions for
these two cases start as
\begin{multline}
\sigma_k^3(n) =\zeta(s)\zeta(s-k)\zeta(s-2k)\zeta(s-3k)
\prod_p
(1+p^{2k-s})^2
(1+p^{k-s})^2
\\ \times
(1-p^{4k-2s})
(1-p^{3k-2s})^3
(1-p^{2k-2s})
(1+p^{6k-3s})^2
(1+p^{5k-3s})^6
\\ \times
(1+p^{4k-3s})^6
(1+p^{3k-3s})^2
(1-p^{8k-4s})^3
(1-p^{7k-4s})^{12}
(1-p^{6k-4s})^{15}
\\ \times
(1-p^{5k-4s})^{12}
(1-p^{4k-4s})^{3}
\cdots,\quad s>1+3k,
\end{multline}
and
\begin{multline}
\sigma_k^4(n) =\zeta(s)\zeta(s-k)\zeta(s-2k)\zeta(s-3k)\zeta(s-4k)
\prod_p
(1+p^{3k-s})^3
\\ \times
(1+p^{2k-s})^5
(1+p^{k-s})^3
(1-p^{6k-2s})^3
(1-p^{5k-2s})^{12}
(1-p^{4k-2s})^{14}
(1-p^{3k-2s})^{12}
\\ \times
(1-p^{2k-2s})^3
(1+p^{9k-3s})^8
(1+p^{8k-3s})^{36}
(1+p^{7k-3s})^{72}
(1+p^{6k-3s})^{88}
(1+p^{5k-3s})^{72}
\\ \times
(1+p^{4k-3s})^{36}
(1+p^{3k-3s})^8
\cdots, \quad s>1+4k.
\end{multline}

\subsubsection{Hybrid Products}
Dirichlet series of mixed products are
\cite[(1.3.3)]{Titchmarsh}\cite{VaidyanTAMS33,BorweinRJ7}
\begin{equation}
\sigma_a(n)\sigma_b(n) \mapsto \frac{\zeta(s)\zeta(s-a)\zeta(s-b)\zeta(s-a-b)}{\zeta(2s-a-b)}
,
\end{equation}
of which
(\ref{eq.Dsigmak2})
is a special case.
An example of this type is
$\sigma_0(n)\sigma_1(n)$ in A064840.

\subsubsection{Dirichlet Convolutions}

$\tau_k(n)$ is the number of ways of expressing $n$ as a product of
$k$ factors.
$\tau_2(n) = \sigma_0(n)$ and iterated convolution with $1$ yield
the ladder of larger $k$
\cite[(1.2.2.)]{Titchmarsh}\cite{VaidyanTAMS33}:
\begin{equation}
\tau_k(n)\mapsto \zeta^k(s).
\label{eq.Dtauk}
\end{equation}
\begin{rem}
$\tau_2$ is A000005,
$\tau_3$ is A007425,
$\tau_4 =\sigma_0(n)\star \sigma_0(n)$ is A007426,
$\tau_5$ is A061200,
  $\tau_6$ is A034695,
  $\tau_7$--$\tau_{11}$ are A111217--A111221,
and  $\tau_{12}$ is A111306.
\end{rem}

\begin{rem}
$\sigma_1^2(n)\star 1$ is A065018.
$\sigma_0^2(n)\star 1$ is A062367.
$\sigma_0^3(n)\star 1$ is A097988.
\end{rem}

\subsection{Powers times Sigma}
\subsubsection{Ordinary product}

Products  with powers have Dirichlet generating functions
derived from
(\ref{eq.shift})
with
(\ref{eq.Dsigmak}) or (\ref{eq.Dsigmak2}):
\begin{equation}
n^t\sigma_k(n) \mapsto \zeta(s-t)\zeta(s-k-t),\quad k\ge 0.
\end{equation}
\begin{equation}
n^t\sigma_k^2(n) \mapsto
\frac{\zeta(s-t)\zeta^2(s-k-t)\zeta(s-2k-t)}{\zeta(2s-2k-2t)}
\end{equation}
\begin{rem}
This concerns $n\sigma_0(n)$ in A038040,
$n^2\sigma_0(n)$ in A034714, and
$n\sigma_1(n)$ in A064987.
\end{rem}

\subsubsection{Dirichlet convolutions}
Convolutions  with powers have Dirichlet generating functions
which are products of (\ref{eq.Dpow}) with
(\ref{eq.Dsigmak}) or (\ref{eq.Dsigmak2}):
\begin{equation}
n^t\star \sigma_k(n) \mapsto \zeta(s)\zeta(s-t)\zeta(s-k),\quad k\ge 0.
\end{equation}
\begin{equation}
n^t\star \sigma_k^2(n) \mapsto
\frac{\zeta(s)\zeta(s-t)\zeta^2(s-k)\zeta(s-2k)}{\zeta(2s-2k)}.
\end{equation}

\begin{rem}
$n\star \sigma_0(n)$ is A007429.
$n\star \sigma_0^2(n)$ is A062369.
$n\star \sigma_1(n)$ is A060640.
$n^2\star \sigma_0(n)$ is A007433.
$n^2\star \sigma_1(n)=n\star \sigma_2(n)$ is A001001.
$n^3\star \sigma_1(n)=n\star \sigma_3(n)$ is A027847.
Multiplication of (\ref{eq.Dsigmak}) with $\zeta(s-t)$ shows \cite[p.\ 285]{Dickson}
\begin{equation}
\sigma_k(n)\star n^t = \sigma_t(n) \star n^k.
\end{equation}
\end{rem}

\subsection{Sums of Odd Divisors}
The master equation for the sum of odd divisors of $n$,
$\sigma_k^{(o)}(n) \equiv \sum_{d|n,d\,\mathrm{odd}}d^k$
is
\begin{equation}
a(2^e)=1;\quad a(p^e) = \left\{
\begin{array}{ll}
\frac{p^{ke+k}-1}{p^k-1}, \quad k>0,\quad p>2. \\
e+1, \quad k=0,\quad p>2. \\
\end{array}
\right.
\end{equation}
The two Bell series for the prime 2 on one hand or any odd prime on the other hand
repeat (\ref{eq.Mzeta}),
\begin{eqnarray}
\sum_{e\ge 0} \frac{1}{2^{es}}&=& \frac{1}{1-2^{-s}}\\
\sum_{e\ge 0} \frac{p^{k(e+1)}-1}{(p^k-1)p^{es}}&=&
\frac{1}{(1-p^{-s})(1-p^{k-s})}
,\quad p>2.
\end{eqnarray}
The Dirichlet series is the interlaced product
\begin{equation}
\sigma_k^{(o)}(n)\mapsto (1-2^{k-s})\zeta(s)\zeta(s-k),\quad k\ge 0.
\end{equation}
\begin{rem}
The OEIS examples are
$k=0$ in A001227,
$k=1$ in A000593,
and $k=2$--$5$ in A050999--A051002.
\end{rem}

\subsection{Euler's Totient}
\subsubsection{Basis function}

The totient $\varphi(n)$ counts numbers $\le n$ and coprime to $n$,
represented by A000010 and its D-inverse A023900.
The master equation is
\begin{equation}
\varphi(p^e) = (p-1)p^{e-1}, \quad e>0.
\label{eq.Mphi}
\end{equation}
The Bell series factorizes in the form \cite[p.\ 111]{Bellman}
\begin{equation}
\varphi(n) = \mu(n)\star n \mapsto \frac{\zeta(s-1)}{\zeta(s)}
\label{eq.Dphi}
\end{equation}

\begin{rem}
The sum of the $k-$th powers of the divisors coprime to $n$, $\varphi_k(n)$,
is generally not multiplicative for $k> 0$. This is easily shown by
finding small indices that violate the defining equation (\ref{eq.Mdef}).
\end{rem}
\begin{rem}
Equation (\ref{eq.Dphi})
has been generalized to define
$\varphi_{k,l}(n) \equiv \mu(n)n^k \star n^l \mapsto \zeta(s-l)/\zeta(s-k)$
\cite{CashwellPJM9}.
$\varphi_{1,2}(n)$ is A002618, $\varphi_{1,3}(n)$ is A000056, $\varphi_{2,3}(n)$ is A053191.
\end{rem}
The square of (\ref{eq.Dphi}) defines
$\varphi^2(n)$ in A127473,
\begin{equation}
a(p^e)= (p-1)^2p^{2e-2},
\label{eq.Mphi2}
\end{equation}
which leads to the Bell series
\begin{equation}
1+\sum_{e\ge 1}\frac{(p-1)^2p^{2e-2}}{p^{es}}
=\frac{1-2p^{1-s}+p^{-s}}{1-p^{2-s}},
\label{eq.Bphisqr}
\end{equation}
and the infinite Euler product
\begin{multline}
\varphi^2(n) \mapsto
\zeta(s-2)
\prod_p
(1-p^{1-s})^{2} (1+p^{-s})(1-p^{2-2s})(1+p^{1-2s})^{2}
\\ \times
(1-p^{3-3s})^{2}
(1+p^{2-3s})^{4} (1-p^{1-3s})^{2} (1-p^{4-4s})^{3} (1+p^{3-4s})^{8} (1-p^{2-4s})^{5}
\\ \times
(1+p^{1-4s})^{2} (1-p^{5-5s})^{6} (1+p^{4-5s})^{16} (1-p^{3-5s})^{16}
\\ \times
(1+p^{2-5s})^{8}
(1-p^{1-5s})^{2} \cdots,\quad s>3.
\label{eq.Dphisqr}
\end{multline}
Cohen defines a multiplicative function $\varphi'$
with a simple master equation build from the product of (\ref{eq.Mmu})
squared and (\ref{eq.Mphi}) \cite{CohenMathZ74}:
\begin{equation}
\varphi'(n) \equiv \mu^2(n)\varphi(n) \mapsto \prod_p (1+p^{1-s}-p^{-s})
.
\end{equation}
This has already been met in (\ref{eq.D2freed}), which 
can be combined into
\begin{equation}
1\star \varphi'(n) = \rad_2(n).
\end{equation}
\begin{rem}
$\varphi'$ is given by the absolute values of A097945.
\end{rem}

\subsubsection{Basic Convolution}
$\varphi\star \varphi$ is A029935.
(\ref{eq.Dmu2}) and (\ref{eq.Dphi}) combine as (A007431, A063659)
\begin{equation}
\mu(n)\star \varphi(n) \mapsto
\zeta(s-1)/\zeta^2(s)
;\quad
\mu^2(n)\star \varphi(n) \mapsto
\zeta(s-1)/\zeta(2s)
.
\end{equation}

\subsubsection{Ordinary product with powers}

The $\varphi(n^t)$ are obtained from (\ref{eq.Dphi}) by the
substitution $e\to et$, so the Bell series is
\begin{equation}
1+\sum_{e\ge 1}\frac{(p-1)p^{et-1}}{p^{es}}
= 
\frac{p^s-p^{t-1}}{p^s-p^t}
=
\frac{1-p^{t-s-1}}{1-p^{t-s}}
\end{equation}
and therefore
\begin{equation}
\varphi(n^t)\mapsto \zeta(s-t)/\zeta(s+1-t).
\end{equation}
Applying (\ref{eq.shift}) establishes the well-known \cite{CohenMathZ74}
\begin{equation}
\varphi(n^t) = n^{t-1} \varphi(n).
\end{equation}
\begin{rem}
This describes
$n\varphi(n)$ (A002618), twice the sum of the integers coprime to $n$ and
not exceeding $n$ \cite{CashwellPJM9}, and
$n^2\varphi(n)$ (A053191).
\end{rem}

\subsubsection{Dirichlet product with powers}

$\varphi\star 1=n$ is obvious from (\ref{eq.Dphi}) \cite{CashwellPJM9}.
Building $\varphi^2\star 1$ we generate A029939\@.
From (\ref{eq.astarbM}) and (\ref{eq.Mphi2}) its master equation ensues,
\begin{equation}
a(p^e)= 1+\sum_{l=1..e} (p-1)^2p^{2l-2} = \frac{p^{2e}(p-1)+2}{p+1}
.
\end{equation}
By construction, the Dirichlet series is (\ref{eq.Dphisqr}) multiplied by $\zeta(s)$.
\begin{rem}
$n\star \varphi(n)$ is A018804 with D-inverse in A101035.
$n^2\star \varphi(n)$ is A069097.
\end{rem}

\subsection{Jordan Functions}

Dirichlet convolution of $n^k$ and $\mu(n)$ defines Jordan functions $J_k$.
The generating functions are an immediate consequence of
(\ref{eq.shift}) and (\ref{eq.Dmu}):
\begin{equation}
n^k \star \mu(n) = J_k(n) \mapsto \zeta(s-k)/\zeta(s).
\label{eq.DJordan}
\end{equation}
\begin{rem}
OEIS representatives are
A000010 ($k=1$),
A007434 ($k=2$) with D-inverse A046970,
A059376 ($k=3$) with D-inverse A063453,
A059377--A059378 ($k=4$--$5$),
and A069091--A069095 ($k=6$--$10$).
\end{rem}

Via (\ref{eq.astarbM}), the master equation for $J_k$ is
\begin{equation}
a(p^e) =p^{k(e-1)}(p^k-1),\quad e>0.
\label{eq.MJordan}
\end{equation}
An immediate consequence of the divisibility properties of the
cyclotomic polynomial $p^k-1$ in this equation is that
$J_d(n)|J_k(n)$ if $d|k$ \cite{BrentMC61,LamAMM103}.

\subsubsection{Products}
A000056 is $nJ_2(n)$.
A115224 is $n^2J_3(n)$.
The convolution products $n^k\star J_k(n)\mapsto \zeta^2(s-k)/\zeta(s)$
generalize Pillai's function \cite{HaukkanenAM76}.

\subsubsection{Dedekind $\psi$}
The Dedekind $\psi$-function is the ratio
\begin{equation}
\psi(n) = J_2(n)/J_1(n) \mapsto \zeta(s)\zeta(s-1)/\zeta(2s)
,
\label{eq.DDedek}
\end{equation}
which can be phrased as 
\begin{equation}
\psi(n) = n\star \xi_t(n)
\end{equation}
with the aid of (\ref{eq.conv}), (\ref{eq.Dpow}) and (\ref{eq.Dsqrfree}).
\begin{rem}
The M\"obius transform $\mu(n)\star \psi(n)$ drops the factor $\zeta(s)$
in (\ref{eq.DDedek}) and is found in A063659.
The Dirichlet series of $n\psi(n)$ (A000082)
and $n^2\psi(n)$ (A033196)
follow from (\ref{eq.shift}).
$\mu(n)\star n\psi(n)$ is A140697.
\end{rem}

\begin{rem}
The $J_k(n)/J_1(n)$ for $k=2$--$17$ are
A001615,
A160889,
A160891,
A160893,
A160895,
A160897,
A160908,
A160953,
A160957,
A160960,
A160972, A161010, A161025, A161139, A161167, and A161213.
\end{rem}

The master equation for $J_k(n)/J_1(n)$ is
a ratio of terms of (\ref{eq.MJordan}):
\begin{equation}
a(p^e) = \frac{p^{k(e-1)}(p^k-1)}{p^{e-1}(p-1)}
=
\frac{p^{(k-1)(e-1)}(p^k-1)}{p-1},\quad e>0,
\end{equation}
with Bell series
\begin{equation}
1+\sum_{e\ge 1}\frac{p^{(k-1)(e-1)}(p^k-1)}{(p-1)p^{es}}
=
\frac{p-1+p^{k-1-s}-p^{-s}}{(p-1)(1-p^{k-1-s})}
=
\frac{1+(\sum_{l=0}^{k-2}p^l)p^{-s}}{1-p^{k-1-s}}
.
\label{eq.DJkJ1}
\end{equation}
At $k=2$ this reduces
to (\ref{eq.DDedek}).
If $k>2$,
(\ref{eq.DJkJ1}) is (\ref{eq.Bsig1nt}) multiplied by $1-p^{-s}$ followed
by the substitution $s\to s+1$;
the Dirichlet series of $J_k(n)/J_1(n)$ are
obtained from prime products like
(\ref{eq.Dsig1n3}) by deleting $\zeta(s)$ and the substitution $s\to s+1$,
to wit
\begin{equation}
1\star [n\frac{J_k(n)}{J_1(n)}] = \sigma_1(n^k),\quad k=1,2,\ldots .
\end{equation}

Multiplicative generalized Dedekind functions $\psi_k=J_{2k}(n)/J_k(n)$
are another generalization 
which---by virtue
of (\ref{eq.MJordan})---have
integer entries governed by
\begin{equation}
a(p^e) = p^{k(e-1)}(p^k+1),\quad e>0.
\end{equation}
The Bell series are
\begin{equation}
1+\sum_{e\ge 1}\frac{p^{ke-k}(p^k+1)}{p^{es}} = \frac{1+p^{-s}}{1-p^{k-s}}
,
\end{equation}
and their product over all primes generates
\begin{equation}
\psi_k\mapsto \zeta(s)\zeta(s-k)/\zeta(2s).
\end{equation}
Mediated by (\ref{eq.Dsqrfree}) and (\ref{eq.Dsigmak}), factorizations of this product
lead to
\begin{eqnarray}
\psi_k(n) &=& n^k\star \mu^2(n); \\
\epsilon_2(n) \star \psi_k(n) &=& \sigma_k(n).
\end{eqnarray}
\begin{rem}
Associated OEIS entries are
$\psi_2(n)$ (A065958),
$\psi_3(n)$ (A065959),
$\psi_4(n)$ (A065960),
$J_6(n)/J_2(n)$ (A194532)
and
$J_8(n)/J_4(n)$ (A194533)
.
\end{rem}

\subsection{Sigma times Totient}
\subsubsection{Ordinary products}

The multiplicative $\sigma_0(n)\varphi(n)$ is represented by A062355\@.
The master equation is the product of (\ref{eq.Mzeta}) and (\ref{eq.Mphi}),
\begin{equation}
a(p^e) = (e+1)(p-1)p^{e-1},\quad e>0.
\end{equation}
The Bell series is
\begin{equation}
1+\sum_{e\ge 1}(e+1)(p-1)p^{e-1}/p^{es}
=
\frac{1-2p^{-s}+p^{1-2s}}{(1-p^{1-s})^2}
.
\label{eq.Bsigma0phi}
\end{equation}
The infinite Euler product is
\begin{multline}
\sigma_0(n)\varphi(n)\mapsto
\zeta^2(s-1)\prod_p
(1-p^{-s})^{2} (1+p^{1-2s})(1-p^{-2s})(1+p^{1-3s})^{2}  \\
\times
(1-p^{-3s})^{2} (1+p^{1-4s})^{4} (1-p^{-4s})^{3} (1-p^{2-5s})^{2} (1+p^{1-5s})^{8}\cdots
,\quad s>2.
\label{eq.Dsigma0phi}
\end{multline}

$\sigma_0(n)\varphi^2(n)$ is A126775 with Bell series
\begin{equation}
1+\sum_{e\ge 1}\frac{(e+1)(p-1)^2p^{2(e-1)}}{p^{es}}
=\frac{1-4p^{1-s}+2p^{3-2s}+2p^{-s}-p^{2-2s}}{(1-p^{2-s})^2}
\end{equation}
and Euler product
\begin{multline}
\zeta_D = \zeta^2(s-2)\prod_p
(1-p^{1-s})^{4} (1+p^{-s})^{2} (1+p^{3-2s})^{2} (1-p^{2-2s})^{7}
\\ \times
(1+p^{1-2s})^{8}
(1-p^{-2s})(1+p^{4-3s})^{8} (1-p^{3-3s})^{28} (1+p^{2-3s})^{34}
\\ \times
(1-p^{1-3s})^{16}
(1+p^{-3s})^{2} (1-p^{6-4s})\cdots,\quad s>3.
\end{multline}

$\sigma_0^2(n)\varphi(n)$ is A110601 with Bell series
\begin{equation}
1+\sum_{e\ge 1}\frac{(e+1)^2(p-1)p^{e-1}}{p^{es}}
=
\frac{p^{1-s}+1-p^{2-3s}+3p^{1-2s}-4p^{-s}}{(1-p^{1-s})^3}
.
\end{equation}
\begin{multline}
\sigma_0^2(n)\varphi(n) \mapsto
\prod_p
(1+p^{1-s})
(1-p^{-s})
(1+p^{1-2s})^{7} (1-p^{-2s})^{6} (1-p^{2-3s})^{8} \\ 
 \times (1+p^{1-3s})^{28} (1-p^{-3s})^{20} (1+p^{3-4s})^{8} (1-p^{2-4s})^{53} (1+p^{1-4s})^{112} \\ 
 \times (1-p^{-4s})^{60} \cdots,\quad s>2.
\end{multline}

The master equation of $\sigma_1(n)\varphi(n)$ (A062354)
is a product of (\ref{eq.Mzeta}) by (\ref{eq.Mphi}),
\begin{equation}
a(p^e) = p^{e-1}(p^{e+1}-1).
\end{equation}
The Bell series is
\begin{equation}
1+\sum_{e\ge 1}
\frac{p^{e-1}(p^{e+1}-1)}{p^{es}}
= \frac{1-p^{1-s}-p^{-s}+p^{2-2s}}{(1-p^{2-s})(1-p^{1-s})},
\end{equation}
and the Euler product
\begin{multline}
\zeta_D(s) = \zeta(s-2)\prod_p
(1-p^{-s})(1+p^{2-2s})(1-p^{1-2s})(1+p^{3-3s})\\ 
 \times (1-p^{1-3s})(1+p^{4-4s})(1+p^{3-4s})(1-p^{1-4s})(1+p^{4-5s})\\ 
 \times (1+p^{3-5s})(1-p^{2-5s})(1-p^{1-5s})\cdots
,\quad s>3.
\end{multline}

\subsubsection{Dirichlet convolutions}
The application of (\ref{eq.shift}) and (\ref{eq.conv})
to (\ref{eq.Dsigmak}) and (\ref{eq.Dphi})
yields \cite[p.\ 293]{Dickson}
\begin{equation}
n^t\varphi(n)\star \sigma_t(n) = n^{1+t}\star 1
\end{equation}
\begin{rem}
Examples of these convolutions are
$\sigma_0(n)\star \varphi(n)$ in A000203,
$\sigma_0^2(n)\star \varphi(n)$ in A060724,
$\sigma_1(n)\star \varphi(n)$ in A038040,
and
$\sigma_2(n)\star \varphi(n)$ in A064987.
\end{rem}

\section{Miscellany}

\subsection{Ramanujan sums}
For our purposes the following definition suffices \cite{Apostol,HaukkanenRJ26}:
\begin{defn} (Ramanujan sum $c_k(n)$)
\begin{equation}
c_k(n) = \sum_{d\mid n,d\mid k} \mu(k/d) d.
\end{equation}
\end{defn}
The associated Dirichlet series are
\cite{Titchmarsh,HardyWright}:
\begin{eqnarray}
c_n(k)\mapsto \frac{\sigma_{1-s}(k)}{\zeta(s)}
,
\\
c_k(n)\mapsto \zeta(s)\sum_{d|k}\mu(k/d)d^{1-s},
\end{eqnarray}
and
\begin{equation}
c_k(n)\tau(n) \mapsto \zeta^2(s)\sum_{\delta|k} \delta^{1-s}\mu(k/\delta)
\prod_{p|\delta}(l+1-lp^{-s})
\end{equation}
where $\delta\equiv \prod p^l$.
\begin{rem}
We find $c_n(1)=\mu(n)$, $c_n(2)$ in A086831, $c_n(3)$ in A085097,
$c_n(4)$ in A085384,
$c_n(5)$ in A085639,
and
$c_n(6)$ in A085906.
$c_1(n)=1$, but
if the role of the argument and index are swapped, the functions
are non-multiplicative in general:
$c_2(n)=-(-1)^n$ and 
$c_3(n)$ in A099837,
$c_4(n)$ in A176742,
and
$c_6(n)$ in A100051.
\end{rem}

\subsection{Unitary Arithmetics}
\subsubsection{Properties}
The unitary convolution
\begin{equation}
(a \oplus b)(n)\equiv \sum_{d\mid n,(d,n/d)=1} a(d)b(n/d)
\end{equation}
shows parallels to the Dirichlet convolution.
Because it preserves the multiplicative
property of its factors \cite{CohenMathZ74,SandorASCM6,TothMP13}
and because its basic associated M\"obius, Sums-of-Divisors and totient functions
are multiplicative, inheritance
similar to Section \ref{sec.core} ensues.
The formula that parallels (\ref{eq.astarbM}) is
\begin{equation}
(a\oplus b)(p^e) = a(1)b(p^e)+a(p^e)b(1),\quad e>0.
\label{eq.astarUbM}
\end{equation}
Cohen defines for example \cite{CohenMathZ74}
\begin{equation}
\sigma'(n) = n\mu^2(n) \oplus 1
.
\end{equation}
Because the master equation of $n\mu^2(n)\mapsto \zeta(s-1)/\zeta(2s-2)$ is
\begin{equation}
a(p^e) = \left\{
\begin{array}{ll}
p^e, & e\le 1 ;\\
0, & e> 1 ,\\
\end{array}
 \right.
\end{equation}
the master equation of $\sigma'$ is constructed from (\ref{eq.astarUbM}) as
\begin{equation}
a(p^e) = \left\{
\begin{array}{ll}
1, & e=0 ;\\
1+p, & e= 1 ;\\
1, & e> 1 .\\
\end{array}
 \right.
\end{equation}
The Bell series is $(1+p^{1-s}-p^{1-2s})/(1-p^{-s})$, which leads to the Dirichlet
series
\begin{multline}
\sigma'(n)\mapsto \zeta(s)\prod_p (1+p^{1-s}-p^{1-2s})
= \zeta(s)\prod_p
(1+p^{1-s})(1-p^{1-2s})
\\ \times
(1+p^{2-3s})
(1-p^{3-4s})(1+p^{4-5s})
(1+p^{3-5s})(1-p^{5-6s})(1-p^{4-6s})\cdots,\quad s>2.
\label{eq.Dsigmaprime}
\end{multline}
\begin{rem}
$\sigma'$ is A092261.
\end{rem}

\subsubsection{Unitary $\mu$}
The $\omega$-analog of (\ref{eq.liouv}) is
the unitary M\"obius function
(A076479)
\cite{CohenMathZ74,SandorASCM6,DaykinQMath17}
\begin{equation}
\mu^*(n) = (-1)^{\omega(n)}
,
\end{equation}
where $\omega(n)$ is the number of distinct prime factors of $n$.
Master equation and Bell series are \cite{CohenAA7}
\begin{equation}
a(p^e) = -1;\quad 1+\sum_{e\ge 1} \frac{a(p^e)}{p^{es}}
= \frac{1-2p^{-s}}{1-p^{-s}}
.
\label{eq.Dmustar}
\end{equation}
The Dirichlet series is $\zeta(s)$ divided by (\ref{eq.eulMc}) at $c=2$,
i. e., $\zeta(s)$ multiplied by the associated Feller-Tornier constant \cite[Tab.\ 6]{MatharArxiv0903}:
\begin{multline}
\mu^*(n)\mapsto
\prod_p
(1-p^{-s}) (1-p^{-2s})(1-p^{-3s})^{2} (1-p^{-4s})^{3} (1-p^{-5s})^{6} \\ 
 \times (1-p^{-6s})^{9} (1-p^{-7s})^{18} (1-p^{-8s})^{30} \cdots,\quad s>1.
\end{multline}
\begin{rem}
The Dirichlet series of Cohen's exponentially odd numbers $\mu_2^*(n)$
is the same at doubled argument $2s$ \cite{CohenMathZ74}.
\end{rem}

\subsubsection{Unitary Sigma}
The unitary $\sigma$-function sums over the divisors $d$ which are
coprime to their complementary divisors $n/d$:
\begin{defn} (Unitary sigma $\sigma^\star$)
\begin{equation}
\sigma_k^*(n)
=
n^k\oplus 1
=\sum_{d|n, (d,n/d)=1}d^k.
\end{equation}
\end{defn}
Applying (\ref{eq.astarUbM}), the master equation for the $k$-power of the divisors is \cite{WallCMB18}
\begin{equation}
a(p^e)=1+p^{ke}.
\label{eq.Musigma}
\end{equation}
The Bell series is
\begin{equation}
1+\sum_{e\ge 1} (1+p^{ek})/p^{es} = \frac{1-p^{k-2s}}{(1-p^{-s})(1-p^{k-s})},
\label{eq.Bsigstar}
\end{equation}
which becomes
\begin{equation}
\sigma_k^*(n) \mapsto \zeta(s)\zeta(s-k)/\zeta(2s-k).
\end{equation}
Multiplication with $\zeta(2s-k)$ generates in view of (\ref{eq.Dz2sk})
and (\ref{eq.Dsigmak})
\begin{equation}
(n^{k/2}\epsilon_2(n))\star \sigma_k^*(n) = \sigma_k(n).
\end{equation}

The sum of the $k$-th power of the odd unitary divisors $\sigma_k^{*(o)}(n)$
is determined by a master equation which 
counts only the first or both of the terms in (\ref{eq.Musigma}) depending on $p$ being even or odd:
\begin{equation}
a(p^e)=\left\{
\begin{array}{ll}
1, & p=2;\\
1+p^{ek}, & p>2.
\end{array}
\right.
\end{equation}
The Bells series is $1/(1-2^{-s})$ for $p=2$ and
(\ref{eq.Bsigstar}) for $p>2$.
In summary
\begin{equation}
\sigma_k^{*(o)}(n)\mapsto  \frac{\zeta(s)\zeta(s-k)(1-2^{k-s})}{\zeta(2s-k)(1-2^{k-2s})}
.
\end{equation}
\begin{rem}
$\sigma_0^*(n)$ is A034444.
$\sigma_1^*(n)$ is
A034448 with D-inverse in A178450.
$\sigma_k^*(n)$ with $k=2$--$8$ are
A034676--A034682.
$\sigma_0^{*(o)}(n)$ is A068068.
$\sigma_1^{*(o)}(n)$ is A192066.
\end{rem}

\subsubsection{Unitary Phi}
The unitary totient
is the unitary convolution of $\mu^*$ and $n$ \cite{CohenMathZ74}:
\begin{defn} (Unitary Totient)
\begin{equation}
\varphi^*(n) = \mu^\star(n) \oplus n
.
\end{equation}
\end{defn}
The master equation is \cite{LalMCom28}
\begin{equation}
a(p^e)= p^e-1
\label{eq.Mutot}
\end{equation}
which sums to
\begin{equation}
1+\sum_{e\ge 1}\frac{p^e-1}{p^{es}} = \frac{1-2p^{-s}+p^{1-2s}}{(1-p^{-s})(1-p^{1-s})}
.
\end{equation}
Comparison of numerator and denominator with (\ref{eq.Bsigma0phi}) shows
that the
Dirichlet series is given by replacing
one of the two $\zeta(s-1)$ in (\ref{eq.Dsigma0phi}) by $\zeta(s)$; this can
be phrased via (\ref{eq.Dphi}) as
\begin{equation}
\varphi^*(n)\star \varphi(n) =
\sigma_0(n)\varphi(n).
\end{equation}
\begin{rem}
$\varphi^*(n)$ is A047994
\end{rem}

The unitary Jordan functions generalize $\varphi^\star(n)$
akin to (\ref{eq.DJordan}) \cite{NageswaraMM70}:
\begin{equation}
J_k^\star(n) = \mu^\star(n)\oplus n^k.
\end{equation}
Via (\ref{eq.astarUbM}), its master equation and Bell series are
\begin{equation}
a(p^e) = p^{ek}-1,
\label{eq.MuJord}
\end{equation}
\begin{equation}
1+\sum_{e \ge 1}\frac{p^{ek}-1}{p^{es}}
=
\frac{1-2p^{-s}+p^{k-2s}}{(1-p^{-s})(1-p^{k-s})}.
\end{equation}
The infinite Euler product becomes
\begin{multline}
J_k^\star(n)
\mapsto
\zeta(s-k)
\prod_p
(1-p^{-s})
(1+p^{k-2s})
(1-p^{-2s})
(1+p^{k-3s})^2
\\ \times
(1+p^{k-4s})^4
(1-p^{-4s})^3
(1-p^{2k-5s})^2
(1+p^{k-5s})^8
(1-p^{-5s})^6
\cdots,\quad s>1+k
.
\end{multline}
(\ref{eq.Mutot}) and (\ref{eq.MuJord}) are related by the substitution $e\to ek$
on the right hand sides, which shows
\begin{equation}
J_k^*(n) = \varphi^*(n^k).
\end{equation}

Unitary analogues of (\ref{eq.Dtauk}) might be created as
\begin{equation}
\tau^*_2(n) = \sigma_0^*(n);\quad \tau_{k+1}^*(n) = \tau_k^*(n)\oplus 1.
\end{equation}
The Bell series is bootstrapped from (\ref{eq.Musigma})
with (\ref{eq.astarUbM}),
\begin{equation}
1+\sum_{e\ge 1}\frac{k}{p^{es}} = \frac{1+(k-1)p^{-s}}{1-p^{-s}}.
\end{equation}
The similarity with (\ref{eq.Bsigma0t}) induces
\begin{equation}
1\star \tau^*_k(n) = \sigma_0(n^k).
\end{equation}

\begin{rem}
$J_1^*(n)$ is A047994.  $J_2^*(n)$ is A191414.
$\tau_3^*(n)$ is A074816.
\end{rem}

\subsection{Higher Order M\"obius}
Apostol's higher order $\mu_k(n)$ generalize (\ref{eq.Mmu}) and are defined
as $\mu_k(n)=0$ if
any prime power $p^{k+1}$ divides $n$, and $\mu_k(n)=(-1)^r$  where $r$ is
the number of maximum prime powers $p^k$ which divide $n$ \cite{Apostol,ApostolPJM32,BegeAUS1}.
The master equation is
\begin{equation}
a(p^e) = \left\{
\begin{array}{ll}
1,& 0\le e<k;\\
-1,& e=k;\\
0,& e>k.\\
\end{array}
\right.
\end{equation}
The Bell series is
\begin{equation}
\sum_{e=0}^{k-1}\frac{1}{p^{es}}-\frac{1}{p^{ks}} = \frac{1-2p^{-ks}+p^{-(k+1)s}}{1-p^{-s}},
\end{equation}
with Dirichlet generating function
\begin{multline}
\mu_k(n)\mapsto \zeta(s)
\prod_p
(1-p^{-ks})^2 (1+p^{-(k+1)s}) (1-p^{-2ks}) (1+p^{-(2k+1)s})^2
\\ \times
(1-p^{-3ks})^2
  (1+p^{-(3k+1)s})^4  (1-p^{-(3k+2)s})^2 (1-p^{-4ks}) ^3(1+p^{-(4k+1)s})^8
\\ \times
(1-p^{-(4k+2)s})^5
(1+p^{-(4k+3)s})^2
(1-p^{-5ks})^6 (1+p^{-(5k+1)s})^{16}
\cdots,\quad s>1.
\end{multline}
\begin{rem}
$\mu_2(n)$--$\mu_4(n)$ are A189021--A189023 in the OEIS \cite{EIS}.
$n\star \mu_2(n)$ is A181549.
\end{rem}

\subsection{Powers
Congruential to Zero}

The number of solutions to $x^t\equiv 0$ (mod $n$)
in the interval $1\le x\le n$
is a multiplicative function with
\cite{BottomleySNJ13}
\begin{equation}
a(p^e) = p^{e-\lceil e/t\rceil} = p^{\lfloor (t-1)e/t\rfloor}.
\label{eq.Mbottom}
\end{equation}
\begin{proof}
It is multiplicative because solutions $x$ for $n$ a product of
prime powers are all products of solutions to the individual prime powers,
and therefore the cardinality of the solutions equals the product
of the cardinality of solutions to the individual prime powers. The
master equation is derived by noting that the solutions
are $x=c p^{\lceil e/t\rceil}$, $c=1,2,\ldots$, with a maximum of
$x=p^e$. The number of solutions equals the maximum $c$, which
is the maximum solution divided by the minimum solution.
\end{proof}
The Bell series is accumulated by splitting
$e=kt+r$ with remainder $0\le r <t$,
and treating $r=0$ and $r\neq 0$ separately:
\begin{multline}
\sum_{e\ge 0} \frac{p^{\lfloor (t-1)(k+r/t)\rfloor}}{p^{es}}
=
\sum_{k\ge 0} \frac{p^{(t-1)k}}{p^{kts}}
+
\sum_{r=1}^{t-1}\sum_{k\ge 0} \frac{p^{(t-1)k+r-1}}{p^{(kt+r)s}}
\\
=
\frac{ 1+ \sum_{r=1}^{t-1} p^{r-1-rs}
}
{1-p^{t-1-st}}
.
\end{multline}
The case $t=2$ is dealt with by plugging $t=2$ into (\ref{eq.Dmaxbtn}).
The Euler product for the case $t=3$ is
\begin{multline}
\prod_p
\frac{1+p^{-s}+p^{1-2s}}{1-p^{2-3s}}
=\zeta(3s-2)
\prod_p
(1+p^{-s})(1+p^{1-2s})(1-p^{1-3s})
\\ \times
(1+p^{1-4s})
(1+p^{2-5s})
(1-p^{1-5s})(1-p^{2-6s})(1+p^{1-6s})\cdots,\quad s>1
.
\end{multline}

\begin{rem}
$t=2$--$4$ are A000188--A000190.
\end{rem}
The associated smallest positive $x$ whose $t$-th power is divisible by $n$
have master equations $a(p^e)= p^{\lceil e/t\rceil}$ and Bell series
\begin{equation}
\sum_{e\ge 0}\frac{p^{\lceil e/t\rceil}}{p^{es}}
=
\sum_{r=0}^{t-1}\sum_{k\ge 0}\frac{p^{k+\lceil r/t\rceil}}{p^{(tk+r)s}}
=
\frac{1+\sum_{r=1}^{t-1}p^{1-rs}}{1-p^{1-ts}}
.
\end{equation}
For $t=2$, the product over primes is 
\begin{equation}
\min_{x>0, x^2\equiv 0 \bmod n} x \mapsto
\prod_p \frac{1+p^{1-s}}{1-p^{1-2s}}
=
\zeta(2s-1)\zeta(s-1)/\zeta(2s-2).
\end{equation}
For $t=3$, a variation of (\ref{eq.Dsigmaprime}) appears:
\begin{multline}
\min_{x>0, x^3\equiv 0 \bmod n} x \mapsto
\prod_p \frac{1+p^{1-s}+p^{1-2s}}{1-p^{1-3s}}
=
\zeta(3s-1)\prod_p (1+p^{1-s}+p^{1-2s})
\\
=
\zeta(3s-1)\prod_p 
(1+p^{1-s})(1+p^{1-2s})(1-p^{2-3s})(1+p^{3-4s})(1-p^{4-5s})\\ 
 \times (1+p^{3-5s})(1+p^{5-6s})(1-p^{4-6s})\cdots, \quad s>2.
\end{multline}
\begin{rem}
These are A019554 and A019555 for $t=2$ and $t=3$,
A053166 for $t=4$, A015052 and A015053 for $t=5$ and $t=6$.
\end{rem}

\bibliographystyle{amsplain}
\bibliography{all}

\end{document}